
\magnification\magstep1
\hsize 15.5truecm
\scrollmode

\overfullrule=0pt
\def\ln{{\rm ln}}\def\Ad{{\rm Ad}}\def\Diff{{\rm Diff}}
\def\diag{{\rm
diag}}\def\tr{{\rm tr}}\def\det{{\rm det}}\def\dimm{{\rm
dim}}
\def\pr{\partial} 

\def\th{\theta}\def\vare{\varepsilon}
\input amssym.def\input amssym.tex

 at 17.28truept
 at 14.4truept
 at 12truept
 at 10.95truept
 at 17.28truept
 at 10truept

\def\title#1{%
\vskip0pt plus.3\vsize\penalty-100%
\vskip0pt plus-.3\vsize\bigskip\vskip\parskip%
\bigbreak\bigbreak\centerline{\bf #1}\bigskip%
}

\def\chapter#1#2{\vfill\eject
\centerline{\bf Chapter #1}
\vskip 6truept%
\centerline{\bf #2}%
\vskip 2 true cm}

\def\section#1#2{%
\def\\{#2}%
\vskip0pt plus.3\vsize\penalty-100%
\vskip0pt plus-.3\vsize\bigskip\vskip\parskip%
\par\noindent{\bf #1\hskip 6truept%
\ifx\empty\\{\relax}\else{\bf #2\smallskip}\fi}}

\def\subsection#1#2{%
\def\\{#2}%
\vskip0pt plus.3\vsize\penalty-20%
\vskip0pt plus-.3\vsize\medskip\vskip\parskip%
\def\TEST{#1}%
\noindent{\ifx\TEST\empty\relax\else\bf #1\hskip 6truept\fi%
\ifx\empty\\{\relax}\else{#2\smallskip}\fi}}

\def\proclaim#1{\medbreak\begingroup\noindent{\bf #1.---}\enspace\sl}

\def\endproclaim{\endgroup\par\medbreak}


\def\comfig#1#2\par{
\medskip
\centerline{\hbox{\hsize=10cm\eightpoint\baselineskip=10pt
\vbox{\noindent #1}}}\par\centerline{ Figure #2}}

\def\figcom#1#2\par{
\medskip
\centerline
{Figure #1}
\par\centerline{\hbox{\hsize=10cm\eightpoint\baselineskip=10pt
\vbox{\noindent #2}}}}
\def\bull{~\vrule height .9ex width .8ex depth -.1ex}


\def\comfig#1#2\par{
\medskip
\centerline{\hbox{\hsize=10cm\eightpoint\baselineskip=10pt
\vbox{\noindent{\sl  #1}}}}\par\centerline{{\bf Figure #2}}}

\def\figcom#1#2\par{
\medskip
\centerline
{{\bf Figure #1}}
\par\centerline{\hbox{\hsize=10cm\eightpoint\baselineskip=10pt
\vbox{\noindent{\sl  #2}}}}}

\def\em{\sl}

\def\\{\hfill\break}

\def\bull{~\vrule height .9ex width .8ex depth -.1ex}

\def\a{\alpha}
\def\b{\beta}
\def\g{\gamma}

\def\d{\delta}

\def\z{\zeta}

\def\la{\lambda}

\def\CC{{\bf C}}

\def\ZZ{{\bf Z}}

\def\cA{{\cal A}}

\def\cD{{\cal D}}
\def\cE{{\cal E}}

\def\cH{{\cal H}}

\def\cL{{\cal L}}
\def\cM{{\cal M}}
\def\cN{{\cal N}}
\def\cO{{\cal O}}

\def\tr{\mathop{\rm tr}\limits}

\def\Hom{\mathop{\rm Hom}\limits}

\def\Ind{\mathop{\rm Ind}}
\def\res{{\rm res}}


\def\la{\lambda}

\def\sqr#1#2{{\vcenter{\hrule height.#2pt%
\hbox{\vrule width.#2pt height#1pt\kern#1pt%
\vrule width.#2pt}%
\hrule height.#2pt}}}

\centerline{\bf Hitchin systems, higher Gaudin operators and $r$-matrices}
\bigskip
\centerline{B. Enriquez and V. Rubtsov}
\medskip
{\bf Abstract.} {\em
We adapt Hitchin's integrable systems to the case of a punctured
curve. In the case of $\CC P^{1}$ and $SL_{n}$-bundles, they are
equivalent to systems studied by Garnier. The corresponding quantum
systems were identified by B. Feigin, E. Frenkel and N. Reshetikhin with
Gaudin systems. We give a formula for the higher Gaudin operators, using
results of R. Goodman and N. Wallach on the center of the enveloping algebras
of affine algebras at the critical level.
Finally we construct a dynamical $r$-matrix for Hitchin systems for a
punctured elliptic curve, and $GL_{n}$-bundles, and
(for $n=2$) the corresponding quantum system.}
\medskip

\section{Introduction.}{}
In [13], N. Hitchin introduced a class of integrable
systems, attached to a complex curve $X$ and a semisimple Lie group $G$.
The problem of quantization of these systems was then connected by
A. Beilinson and V. Drinfel'd to the Langlands program. This quantization
makes use of differential operators on the moduli space of $G$-bundles
on $X$, obtained from the action of the center of the local completion
of the enveloping algebra of a Kac-Moody algebra, at the critical
level.

This program can also be carried out in the case of curves
with marked points. In the particular case of the punctured $\CC P^{1}$,
the action of the center of the enveloping algebra was studied by B.
Feigin, E. Frenkel and N. Reshetikhin in [6]; they obtained an
integrable system whose first operators are identical to Gaudin's
operators ([9]).

In this paper, we consider the question of expressing the action of
higher central elements. We first note, that the Adler-Kostant-Symes (AKS)
scheme, which allows to write families of commuting operators ([2],
[14], [21]), can be applied in the present situation, and then show that
the higher Hamiltonians obtained in [6], coincide with those. So our
problem turns out to be equivalent to expressing higher central elements
in the enveloping algebras at critical level, a problem which was
solved by several authors ([10], [12]). Here we show how the results
of [10] can be used to derive a simple expression of higher Gaudin
Hamiltonians.

We then turn to the case of punctured elliptic curves. We show that
the integrability of Hitchin's system can be deduced
from an $r$-matrix argument. Here $r$-matrix relations contain
additional terms, due to an invariance under the Cartan algebra action.
The $r$-matrix depends on the moduli parameters, so it reminds
dynamical $r$-matrices. In the case of one puncture, our $L$-operator
and $r$-matrix seem connected with those considered respectively
by I. Krichever and A. Gorsky and N. Nekrasov in [15] and [11], and
H. Braden, T. Suzuki and E. Sklyanin [5],
[19].  It is also analogous to the $r$-matrix appearing in
the work of G. Felder and C. Wieczerkowski on the
Knizhnik-Zamolodchikov-Bernard equation on elliptic curves ([7]).
We give the form of the first Hamiltonians in this case; one of them
contains a Weierstass potential, and so is analogous to the
Calogero-Moser Hamiltonian. We compute the corresponding quantum
Hamiltonians, in the case $G=GL_{2}$.

We would like to thank V. Drinfel'd, B. Feigin, G. Felder, E. Frenkel,
A. Gorsky, N. Nekrasov,
A. Reyman, and A. Stoyanovsky for discussions connected with the subject
of this work. We are thankful to A. Beilinson and V. Drinfel'd for
sending us their paper [4]. The work of V.R. was supported by the CNRS,
and partially by grant RFFI 95-01-01101; he is grateful to the Centre de
Math\'ematiques de l'\'Ecole Polytechnique, where this work was done,
for hospitality.

\section{1.}{Hitchin and Beilinson-Drinfeld systems in the case of a
general punctured curve.}

\medskip\noindent
\it 1.1. Hitchin systems. \rm
Let $\overline X$ be a smooth compact complex curve, $x_{i}$
($i=1,\cdots,N$) be distinct points on $\overline X$. Set $X=\overline
X-\{x_{i}\}$.
Let $G$ be a reductive complex Lie group, $B\subset G$ and $T\subset B$
Borel and Cartan subgroups of $G$; let ${\bf g}$, ${\bf b}$ and ${\bf
t}$ be their Lie
algebras. Let $\cM_{G}(X)$ be the moduli space
of principal $G$-bundles on $\overline X$ with choices of a $B$-orbit in each
fibre over $x_{i}$. Let us identify ${\bf g}$ with
its dual, using a non-degenerate invariant form $\langle,
\rangle_{{\bf g}}$. Let $P
\in\cM_{G}(X)$, then $T^{*}_{P}\cM_{G}(X)$ is formed of the
$\xi\in H^{0}(\overline X,\Omega_{\overline X}(\sum_{i=1}^{N}(x_{i}))\otimes
{\bf g}_{P})$, such that
$\xi$ has the expansion near $x_{i}$, $\xi=
{1\over{u_{i}}}\xi_{i}+ {\rm regular}$, and $\xi_{i}\in
({\bf b}_{x_{i}})^{\perp}$; ${\bf b}_{x_{i}}$
is the subspace of the fibre of ${\bf g}_{P}$
at $x_{i}$, corresponding to the $B$-orbit at $P_{x_{i}}$, $u_{i}$ is a
local coordinate at $x_{i}$.
The Hitchin map
$$
H : T^{*}\cM_{G}(X)\to \cH_{X}=
\bigoplus_{i=1}^{r}H^{0}\big(\overline X,\Omega_{\overline X}^{d_{i}}
((d_{i}-1)\sum_{l=1}^{N}(x_{l}))\big),
$$
is then
defined by $(H(P,\xi))_{l}= P_{d_{l}}(\xi)$; $r$ is the rank of $G$,
$d_{l}$ ($1\le l\le r$)
are the characteristic degrees of ${\bf g}$ and $P_{d_{l}}$ corresponding
invariant polynomials.
The generic fiber of the
natural projection $\cM_{G}(X)\to\cM_{G}(\overline X)$
is $(G/B)^{N}$ if genus$(\overline X)>1$, the generic bundle having no
automorphisms; on the other hand, we have for genus$(\overline X)>0$, $\dimm
\cH_{X} =\dimm\cH_{\overline X}+\sum_{l=1}^{r}(d_{l}-1)N
=\dimm\cH_{X}+N(\dimm B-r)$.
If genus$(\overline X)=1$, an open subset of $\cM_{G}(X)$ is identified
with $T/W$ ($W$
is the Weyl group of $G$) if $N=0$, and with
$T\rtimes W\setminus [T\times(G/B)^{N}]$ for $N\ge 1$
(only $W$ acts in the first factor, and $T\rtimes W$ acts diagonally on
$(G/B)^{N}$); on the other
hand, $\dimm\cH_{X}=\sum_{i=1}^{r}(d_{i}-1)N$ if $N\ge 1$, and $r$
if $N=0$.
If genus$(\overline X)=0$ and $N\ge 3$, an open subset of
$\cM_{G}(X)$ is identified
with $G\setminus(G/B)^{N}$, whereas $\dimm \cH_{X}=\sum_{l=1}^{r}
[(d_{l}-1)(N-2)-1]$. The cases $N\le 2$ give trivial moduli spaces and
$\cH_{X}$. So in all cases
$$
\dimm\cM_{G}(X)= \dimm\cH_{X}.
$$

We can see as in [13] that the functions on $T^{*}\cM_{G}(X)$,
defined by the coordinates of $H$, Poisson commute. Moreover, for
$G=GL_{n}(\CC)$ we can consider the spectral cover of
$\overline X$, defined as
$\{(x,\lambda)| \lambda^{n}+ \sum_{l\ge 1}H_{i}\lambda^{n-l} = 0\}$,
for $(H_{l})\in \cH$ fixed, in the total space of
$\Omega_{\overline X}(\sum_{i=1}^{N}(x_{i}))$; it has ramification of order
$n$ at the
points $x_{i}$, in the generic situation. It is possible to build a line
bundle over the spectral cover, and to
study the integrability of the system as in [13].

\medskip\noindent
\it 1.2. Beilinson-Drinfeld systems. \rm
To quantize the Hitchin systems, Beilinson and Drinfeld ([4]) define
$\dimm\cM_{G}(X)$ commuting differential operators on $\cM_{G}(X)$, with
symbols the coordinates of the map $H$ (here we assume no marked
points). They are constructed as follows:
a base point $x$ on $X$ being fixed, $\cM_{G}(X)$ is identified with
$G(\cO_{x})\setminus G(k_{x})/G(A)$ (Siegel-Weil); $\cO_{x}$ and $k_{x}$
are respectively the local ring and field at $x$, and $A=H^{0}(X-\{x\},
\cO_{X})$. Then the center $Z(U_{-h^{\vee}}({\bf g}(k_{x}))_{loc})$
of $U_{-h^{\vee}}({\bf g}(k_{x}))_{loc}$ (local
completion of the enveloping algebra of the critical level extension of
${\bf g}(k_{x})$) acts by differential operators on the line bundle
$(\det)^{-h^{\vee}}$ over $\cM_{G}(X)$. This procedure can easily be
extended to the punctured case: remark that
$\cM_{G}(X)=
G(\cO_{x})\setminus G(k_{x})/\Gamma$, where $\Gamma\subset G(A)$  is
formed of the regular maps from $\overline X-\{x\}$ to $G$, taking values in
$B$ at points $x_{i}$. Let $(\la_{i})_{1\le i\le N}$ be a system of
dominant weights of $G$. We define a line bundle $\cL_{(\la_{i})}$ on
$\cM_{G}(X)$ as follows: $(\la_{i})$ defines a character of $\Gamma$ (by
the maps $\Gamma\to B^{N}\to T^{N}$) and so a line bundle
$\cL'_{(\la_{i})}$ on $G(k_{x})/\Gamma$, then $\cL'_{(\la_{i})}\otimes
(\det)^{-h^{\vee}}$ has a natural action of $G(\cO_{x})$;
$\cL_{(\la_{i})}$ is then the quotient
bundle. The center of $U_{-h^{\vee}}({\bf g}(k_{x}))_{loc}$ then acts on this
bundle by differential operators as before.

\section{2.}{Hitchin systems in the rational case.}

In this section and the following, we set $\overline X=\CC P^{1}$, and
denote by $z_{i}$ the coordinate of the marked point $x_{i}$
($i=1,...,N$); we assume that no $x_{i}$ coincides with $\infty$.

We will express the corresponding Hitchin systems, and recall an $r$-matrix
result of Semenov about them.

An open subset of $\cM_{G}(X)$ is formed by parabolic structures on the
trivial bundle; this subset, that we call $\cM_{G}^{(0)}(X)
$ is isomorphic to $G\setminus (G/B)^{N}$ [the left action of $G$ is
diagonal].

Recall the Springer resolution $T^{*}(G/B)\to\cN$, $\cN$ the nilpotent
cone of ${\bf g}$ ([20]). Then we construct, by reduction, the resolution
$$
T^{*}[G\setminus
(G/B)^{N}]\to\{(\eta^{(i)})\in\cN^{N}|\sum_{i=1}^{N}\eta^{(i)}=0\}/G
$$
(the action of $G$ on the last term is by conjugation).
When the $\eta^{(i)}$ are regular, the parabolic structure
corresponding to $(\eta^{(i)})_{i=1,\cdots,N}$ is
$(g_{i}B)_{i=1,\cdots,N}$, where $g_{i}$ are elements of $G$
conjugating $\eta^{(i)}$ to elements of ${\bf b}\subset {\bf g}$. The
$1$-form $\xi$ is then
$$
\xi=\sum_{i=1}^{N}{\eta^{(i)}\over{z-z_{i}}}dz.\leqno(1)
$$
The Poisson structure on $T^{*}\cM_{G}^{(0)}(X)$
corresponds, in terms of the $(\eta^{(i)})$, to the product of
Kostant-Kirillov structures on each $\cN$. In tensor notation: $\{
\eta^{(i)}\otimes_{,}\eta^{(j)}\}=\delta_{ij}[P,1\otimes\eta^{(j)}]
=-\delta_{ij}[P,\eta^{(i)}\otimes 1]$, $P$ the permutation operator.
We deduce from this:
$$
\{\eta(z)\otimes_{,}\eta(w)\}=[{P\over{z-w}},\eta(z)\otimes
1+1\otimes \eta(w)],\leqno(2)
$$
where $\eta(z)=\sum_{i=1}^{N}{\eta^{(i)}\over{z-z_{i}}}$. So we have:

\proclaim{Proposition 2.1} ([18]) Let us endow $\cN^{N}$ with the
product of Kostant-Kirillov structures on each factor. Then the mapping
$\cN^{N}\to {\bf g}[[z^{-1}]]$, $(\eta^{(i)})_{1\le i\le N}
\mapsto \eta(z)=\sum_{i=1}^{N}{\eta^{(i)}\over{z-z_{i}}}$, is
Poisson. \endproclaim

We deduce from this that the coefficients of the forms
$P_{d_{i}}(\eta(z))$ are in involution, because all the central
functions on ${\bf g}[[z^{-1}]]$ are in involution. (This gives
another proof of
involutivity of Hitchin's Hamiltonians.) Let
us give now the expression of the corresponding flows:

\proclaim{Proposition 2.2} Decompose $P_{d_{l}}(\eta(z))$ as
$\sum_{a_{1}+\cdots+a_{N}=d_{l}-1} {H_{d_{l},(a_{i})}\over{
\prod_{i=1}^{N}(z-z_{i})^{a_{i}}
}}(dz)^{d_{l}}$, and denote by $\pr_{d_{l},(a_{i})}$ the flow generated
by $H_{d_{l},(a_{i})}$. Then we have the identity of rational functions
in $z$
$$
\sum_{a_{1}+\cdots+a_{N}=d_{l}-1}
{\pr_{d_{l},(a_{j})}(\eta^{(i)})
\over{\prod_{j=1}^{N}(z-z_{j})^{a_{j}}}}
=
[P'_{d_{l}}(\eta(z)),{\eta^{(i)}\over{z-z_{i}}}].
\leqno(3)
$$
For ${\bf g}=sl_{n}(\CC)$, the r.h.s. is
$[d_{l}(\sum_{j=1}^{N}{\eta^{(j)}\over{z-z_{j}}})^{d_{l}-1},\eta^{(i)}]$.
For ${\bf g}$ arbitrary, the flows corresponding to $d_{1}=2$ are
$$
\pr_{i}\eta^{(j)}=-{[\eta^{(i)},\eta^{(j)}]\over{z_{i}-z_{j}}}
{\rm \ \ for \ \ }j\ne i, {\rm \ \ and \ \ }
\pr_{i}\eta^{(i)}=\sum_{j\ne i}{[\eta^{(i)},\eta^{(j)}]
\over{z_{i}-z_{j}}}.
\leqno(4)
$$
\endproclaim

We note that in the case $g=sl_{n}(\CC)$, the flows $\pr_{i}$ already appeared
in [8] (we thank J. Harnad for pointing out this reference to
us). Their integration was studied by many authors (cf. e.g. [1], [3]).

\section{3.}{Gaudin systems.}

\medskip\noindent
\it 1. The moduli stack in the rational case. \rm

Let again ${\bf g}$ be an arbitrary reductive complex Lie algebra, $G$
be the adjoint group. Let $\Delta$ be the set of the roots of ${\bf g}$
w.r.t.  ${\bf t}$, ${\bf g}={\bf t}\oplus
\bigoplus_{\a\in\Delta}{\bf g}_{\a}$ the associated decomposition of
${\bf g}$, $R$ be
the root lattice of ${\bf g}$, $W$ the Weyl group of ${\bf g}$. Classes
of principal $G$-bundles on $\overline X=\CC P^{1}$ are parametrized by
$\Hom(R,\ZZ)/W$; to $\varpi\in\Hom(R,\ZZ)$ we associate the Lie algebra
bundle on $\CC P^{1}$
$$
{\bf g}(\varpi)={\bf t}
\oplus\bigoplus_{\a\in\Delta}{\bf g}_{\a}(\varpi(\a)\infty),\leqno(5)
$$
and the associated $G$-bundle $G(\varpi)=\Ad {\bf g}(\varpi)$. Its
automorphism group is a subgroup of $G(\CC[z])$,
$$
P_{\varpi}=L_{\varpi}U_{\varpi},
\quad U_{\varpi}=\prod_{\a\in\Delta}N_{\a}(H^{0}(\varpi(\a)\infty)),
\leqno(6)
$$
$L_{\varpi}$ the subgroup of ${\bf g}$ with Lie algebra
${\bf l}_{\varpi}={\bf t}\oplus\bigoplus_{\varpi(\a)=0}{\bf g}_{\a}$.
The moduli space of $G$-bundles on $\CC P^{1}-\{z_{i}\}$ is
$$
\cM_{G}(X)=\prod_{[\varpi]\in\Hom(R,\ZZ)/W}
\cM_{G}^{\varpi}(X),
\leqno(7)
$$
where $\cM_{G}^{\varpi}(X)$ is identified with
$P_{\varpi}\setminus (G/B)^{N}$, where the action of $G(\CC[z])$ is the
composition of the morphism $G(\CC[z])\to G^{N}$, $g(z)\mapsto
(g(z_{i}))_{i}$, and the left translation. Let $(\la_{i})_{i}$ be
integral dominant weights of $G$, $\cL_{\la_{i}}$ be the associated line
bundles on $G/B$; $\boxtimes_{i=1}^{N}\cL_{\la_{i}}$ is a
$G^{N}$-equivariant bundle on $(G/B)^{N}$, so it is
$P_{\varpi}$-equivariant; let $\cL_{(\la_{i})}$ be the quotient bundle
on $\cM_{G}^{\varpi}(X)$.

\medskip\noindent
\it 2. The FFR scheme. \rm

The procedure of sect. 1.2 was applied in [6] to the case of the
punctured $\CC P^{1}$. Let us set some notations.
Let for each $i$, $k_{i}$ and $\cO_{i}$ be the local field and ring at
$z_{i}$; let $\tilde {\bf g}$ the
central extension of $\oplus_{i}{\bf g}(k_{i})$ by the cocycle
$c((a_{i}),(b_{i}))=\sum_{i=1}^{N}\res_{x_{i}}
\langle a_{i}, db_{i}\rangle_{{\bf g}} K$,
with values in the abelian algebra $\CC K$. Let $\tilde {\bf g}_{+}$
be the preimage of ${\bf g}(\oplus_{i}\cO_{i})$ in $\tilde {\bf g}$;
$\tilde {\bf g}_{+}$ is then
isomorphic to ${\bf g}(\oplus_{i}\cO_{i})\oplus
\CC K$. Let for $\la$ integral dominant weight of ${\bf g}$,
$V_{\la}$ be the corresponding irreducible representation of ${\bf g}$;
and let for $k\in \CC$, and $\la_{1},...,\la_{N}$ integral dominant
weights of ${\bf g}$, $V_{(\la_{i})}^{k}$ be the representation of
$\tilde {\bf g}_{+}$ in $V_{\la_{1}}\otimes...\otimes V_{\la_{N}}$,
where elements of ${\bf g}(\oplus_{i}\cO_{i})$ act as their
images in ${\bf g}^{\oplus N}$, and $K$ by $k$.
Let $\bar {\bf g}_{(z_{i})}$ be the Lie algebra of
regular maps from $X$ to ${\bf g}$; choose and denote the same way a
lifting of this algebra to $\tilde {\bf g}$.

Let $\hat {\bf g}$ be the universal central extension of ${\bf g}((u))$,
and ${\bf V}_{0}^{-h^{\vee}}$ be the critical level vacuum module over
it. Central
fields $T(\zeta)\in Z(U_{-h^{\vee}}(\hat{\bf g})_{loc})[[\zeta^{\pm1}]]$
are in
correspondance with imaginary weight singular vectors $\sum
I(-l)J(-k)...v_{0}\in{\bf V}_{0}^{-h^{\vee}}$, $I,J,...\in {\bf g}$.
Following [6], the action of $T(\zeta)$ on $H^{0}(\cM_{G}(X),
\cL_{(\la_{i})})$ can be described as follows. We have an identification
$$
H^{0}(\cM_{G}(X), \cL_{(\la_{i})})=\bar
H_{(\la_{i})}^{-h^{\vee}}=\{\mu\in
({\bf V}_{(\la_{i})}^{-h^{\vee}})^{*}|\mu {\rm \ is \
}\bar {\bf g}_{(z_{i})}{\rm -invariant}\},\leqno(8)
$$
where for any $k$, ${\bf V}^{k}_{(\la_{i})}$ is the induced
representation $\Ind_{\tilde {\bf g}_{+}}^{\tilde {\bf g}}
V^{k}_{(\la_{i})}$.
According to [6], 3, lemma 1,
$$
\bar H_{(\la_{i})}^{-h^{\vee}}\simeq (V_{\la_{1}}\otimes\cdots\otimes
V_{\la_{N}})^{*}.
\leqno(9)
$$
The field $T$ corresponds to an imaginary weight singular vector $\sum
I(-l)J(-k)...v_{0}\in{\bf V}_{0}^{-h^{\vee}}$, $I,J,...\in {\bf g}$. Due
to the ``swapping''
identity (3.1) of [6], the action of this singular vector on
$(V_{\la_{1}}\otimes\cdots\otimes V_{\la_{N}})^{*}$ is
$$
\sum{1\over{(l-1)!}}\pr^{l-1}I(u){1\over{(k-1)!}}\pr^{l-1}J(u)...,
\leqno(10)
$$
where $I(u)=\sum_{i=1}^{N}{{I^{(i)}}\over{u-z_{i}}}$.

For example, the operators corresponding to the degree two Casimir
element are the Gaudin Hamiltonians $H_{2,i}$, such that the combination
$H_{2}(\zeta)=\sum_{i=1}^{N}{H_{2,i}\over{\zeta-z_{i}}}$ satisfies
$$
H_{2}(\zeta)=\sum_{i} e_{i}(z)e_{i}(z),
\leqno(11)
$$
with $(e_{i})$ an orthomormal basis of ${\bf g}$.

\medskip\noindent
\it 3. The AKS scheme. \rm

On the other hand, the expression (1) gives a realization of the
Lie algebra $u^{-1}{\bf g}[[u^{-1}]]$. More precisely, we have a Lie
algebra morphism $\pi:u^{-1}{\bf g}[[u^{-1}]]\to {\bf g}^{\oplus N}$,
defined by
$\pi(Iu^{-k})=\sum_{i=1}^{N}I^{(i)}z_{i}^{k-1}$. Let us show how the AKS
scheme enables us to construct a commuting family in
$U(u^{-1}{\bf g}[[u^{-1}]])$. Let us decompose the central extension
$\CC K\to
\hat {\bf g}\to {\bf g}((u))$ as $\hat {\bf g}={\bf a}\oplus {\bf b}$,
${\bf a}=\sigma(u^{-1}{\bf g}[[u^{-1}]])$
and ${\bf b}=\a^{-1}({\bf g}[u])$, $\a$ being the projection and
$\sigma$ being a
section of $u^{-1}{\bf g}[[u^{-1}]]$ to $\hat {\bf g}$. Then, $U\hat
{\bf g}=U{\bf a}\oplus (U\hat {\bf g}){\bf b}$. We have then an algebra
morphism $Z(U\hat {\bf g})\to U{\bf a}$,
given by the projection to the first factor, whose image is a commuting
family in $U{\bf a}$. Let us specialize this construction to the critial
level. We have then a sequence of morphisms
$$
Z(U_{-h^{\vee}}\hat {\bf g})\to U(u^{-1}{\bf g}[[u^{-1}]])\to
(U{\bf g})^{\otimes N},
$$
the last one being given by $\pi$. This gives a family of commuting
differential operators on $(G/B)^{N}$.

\noindent
\it Remark. \rm  According to [17], Gaudin systems can be obtained from
quantum tops systems by a reduction procedure, which explains that the
AKS scheme can be applied to them.

\medskip\noindent
\it 4. Coincidence of the AKS and FFR systems. \rm

To see that these operators are the same as those obtained by the
previous construction, let us work out the AKS scheme more
explicitly. The central field $T(u)$ associated to $\sum I(-l)J(-k)...$,
is the normally ordered product
$$
\sum{1\over{(l-1)!}}{1\over{(k-1)!}}...(\partial^{l-1}\bar
I(u)(\partial^{k-1}\bar J(u)...()))\leqno(12)
$$
(where the parenthesis stand for the normal ordering operation); here
$$
\bar I(u)=\sum_{n\in\ZZ}I(-n-1)u^{n}=I_{+}(u)+I_{-}(u),\leqno(13)
$$
$I_{+}(u)=\sum_{n\ge 0}I(-n-1)u^{n}$. The transform of this expression
by the AKS procedure will be, due to the conventions
$(AB)(u)=(A_{+}B+BA_{-})(u)$,
$$
\sum{1\over{(l-1)!}}{1\over{(k-1)!}}...(\partial^{l-1}
I_{+}(u)(\partial^{k-1}J_{+}(u)...())).
\leqno(14)
$$
But the image by $\pi$ of $I_{+}(u)$ is
$I(u)=\sum_{i=1}^{N}{{I^{(i)}}\over{u-z_{i}}}$; this shows

\proclaim{Proposition 3.1}
The expressions (10) and (14) for AKS and FFR Hamiltonians, coincide.
\endproclaim

\medskip\noindent
\it 5. Application: expression of the higher Gaudin operators in the
$sl_{n}$ case. \rm

The following can then be deduced from [10], using the Newton
identities.

\proclaim{Proposition 3.2} Let $(s_{p})_{p\ge 0}$ be the sequence of
polynomials in $n$, defined by $s_{1}=0$, $s_{2}=n/2$, $s_{3}=-2n/3$,
and for $p\ge 2$,
$$
(n-p)s_{p}-2(p+1)s_{p+1}-(p+2)s_{p+2}=0;\leqno(15)
$$
let $\la_{1},\cdots,\la_{n}$ be the solutions of the equation
$$
\la_{n}-s_{1}(n)\la^{n-1}+s_{2}(n)\la^{n-2}-s_{3}(n)\la^{n-3}\cdots=0,
\leqno(16)
$$
and let $H=\diag(\la_{1},\cdots,\la_{n})$. Let us set,
${}^{k}H=\Ad(k)H$, for $k\in K=SU(n,\CC)$; and let $dk$ be a Haar
measure on $K$. Then the higher Gaudin
Hamiltonians are the operators $H_{l,a_{i}}$, ($\sum_{i=1}^{N}a_{i}=l-1$),
defined for each $l=2,...,N$ by
$$
\sum_{a_{1}+...+a_{N}=l-1}
{H_{l,a_{i}}\over{\prod_{i=1}^{N}(\zeta-z_{i})^{a_{i}}}}=\int_{K}\bigg(
{{({}^{k}H)^{(i)}}\over{\zeta-z_{i}}}\bigg)^{l}dk.\leqno(17)
$$
\endproclaim

\section{4.}{An $r$-matrix for the case of punctured elliptic curves.}

Let us turn now to the case where $\overline X$ is an elliptic curve
$\CC^{\times}/q^{\ZZ}$ ($q<1$). We denote by $z_{i}$ ($i=1,\cdots,N$)
the coordinates of the marked points. We fix from now on,
$G=GL_{n}(\CC)$.
Consider the open subset $\cM^{(0)}_{G}(\overline X)$ of
$\cM_{G}(\overline X)$, formed of the space of bundles on $\overline X$,
direct
sums of line bundles of degree $0$. These bundles are parametrized by
the symmetric product $\overline X^{(n)}$; to
$(t_{1},\cdots,t_{n})\in (\CC^{\times})^{n}$, we associate the bundle
$\cE_{(t_{\a})}=\CC^{\times}\times\CC^{n}/\{(z,\xi)\sim
(qz,\diag(t_{\a})\xi)\}$ over $X$; changing $(t_{\a})$ into
$(q^{a_{\a}}t_{\a})$ (with the $a_{\a}$ integers)
gives an isomorphic bundle, the isomorphism being $(z,\xi)\mapsto
(z,\diag(z^{a_{\a}})\xi)$.

Now, the preimage in $\cM_{G}(X)$ of this open subset can be
identified with
$$
T\rtimes S_{n}\setminus (\CC^{\times})^{n}\times
(G/B)^{N}/[(t_{\a},g_{i}B)
\sim (q^{a_{\a}}t_{\a},\diag(z_{i}^{a_{\a}})g_{i}B)],
$$
$T\rtimes S_{n}$ acting diagonally on $(G/B)^{N}$, and by permutations on
$(\CC^{\times})^{n}$. We denote it by $\cM^{(0)}_{G}(X)$.
The cotangent to $\cM^{(0)}_{G}(X)$ is now the quotient by
$S_{n}$ of the
reduction by $T$ of $T^{*}((\CC^{\times})^{n}\times (G/B)^{N})$.
The Springer resolution gives now a mapping from
$T^{*}\cM^{(0)}_{G}(X)$ to
$$
\eqalign{
T\rtimes
S_{n}\setminus\{(p_{\a},t_{\a},\eta_{i})
\in\CC^{n}\times(\CC^{\times})^{n}
 & \times\cN^{N} |
(\sum_{i=1}^{N}\eta_{i})_{t}=0\}/
\cr &
/ \{(p_{\a},t_{\a},\eta_{i})\sim(p_{\a},
q^{a_{\a}}t_{\a},  \Ad(\diag(z_{i}^{a_{\a}}))\eta_{i})\},
\cr}$$
bijective over the open subset, defined by the condition that each
$\eta_{i}$ be regular. This map is
compatible with the Poisson bracket given by
the product of $\{p_{\a},t_{\b}\}=\d_{\a\b}t_{\b}$, and Kostant-Kirillov
on each copy of $\cN$.

The corresponding $1$-form $\xi\in
H^{0}(X,gl(\cE_{(t_{\a})})(-\sum_{i=1}^{N}(z_{i})))$ can be seen as a
$1$-form $\tilde{\xi}$ on $\CC^{\times}$ with values in $gl_{n}(\CC)$,
with simple poles at $z_{i}q^{\ZZ}$, and such that
$\tilde{\xi}(qz)=\Ad(t_{\a})\tilde{\xi}(z)$; it is given by
$\tilde{\xi}(z)=\bar{\xi}(z){{dz}\over z}$, with
$$
\bar{\xi}(z)_{\a\b}=\sum_{i=1}^{N}\eta_{\a\b}^{(i)}
{{\th(t_{\a}^{-1}t_{\b}z z_{i}^{-1})}\over {\th(t_{\a}^{-1}t_{\b})\th(z
z_{i}^{-1})}} {\rm \ if } \a\ne \b,
\bar{\xi}(z)_{i}={1\over\th'(1)}p_{\a}+\sum_{i=1}^{N}{1\over\th'(1)}
{\dot\th\over\th}(zz_{i}^{-1})\eta_{i}^{i}.
\leqno(18)
$$
Here $\th(z)=\prod_{i\ge 0}(1-q^{i}z)\prod_{i\ge 1}(1-q^{i}z^{-1})$;
$\th$ has the properties $\th(qz)=-z^{-1}\th(z)$,
$\th(z^{-1})=-z^{-1}\th(qz)$; we denote $\dot\th(z)=z{d\th\over dz}(z)$.

We will show that the commutativity of the coordinates of the
$\tr\tilde{\xi}(z)^{k}$ in a basis of the space of $k$-forms
on $X-\{z_{i}\}$ can be deduced, as in prop. 2.1, from an $r$-matrix
argument:

\proclaim{Proposition 4.1}
Let $r(z,w,t_{\a})$ and $\rho(z,w,t_{\a})$ be the matrices acting on
$\CC^{n}\otimes\CC^{n}$, with elements
$$
r(z,w,t_{\a})_{\a\b}^{\g\d}=\bigg(
-{\th(t_{\a}^{-1}t_{\b}zw^{-1})\over {\th(t_{\a}^{-1}t_{\b})\th(z
w^{-1})}}\d_{\a\b}^{\d\g}
+{1\over\th'(1)}{\dot\th\over\th}(zw^{-1})\d_{\a\b}^{\g\d}
\bigg)(1-\d_{\a\b})\leqno(19)
$$
and
$$
\rho(z,w,t_{\a})_{\a\b}^{\g\d}={1\over\th'(1)}
{\th(t_{\a}^{-1}t_{\b}zw^{-1})\over {\th(t_{\a}^{-1}t_{\b})\th(z
w^{-1})}}
\bigg[
{\dot\th\over\th}(t_{\a}t_{\b}^{-1})+{\dot\th\over\th}(t_{\a}^{-1}t_{\b}
zw^{-1}) \bigg]
\d_{\a\b}^{\d\g}(1-\d_{\a\b});\leqno(20)
$$
let $\bar{\xi}(z)$ be given by formula (18); let us endow the system of
variables $(p_{\a},t_{\a},\eta_{i})$  with the Poisson brackets, product
of  $\{p_{\a},t_{\b}\}=\d_{\a\b}t_{\b}$ and Kostant-Kirillov on each copy of
$\cN$; then we have
$$
\eqalign{
\{\bar\xi(z,t_{\a})\otimes_{,}\bar\xi(w,t_{\a})\}=
[r(z,w,t_{\a}), & \bar\xi(z,t_{\a})\otimes 1+1\otimes \bar\xi (w,t_{\a})]
\cr & +\rho(z,w,t_{\a})[(\sum_{i=1}^{N}\eta_{i})_{t}\otimes 1
-1\otimes(\sum_{i=1}^{N}\eta_{i})_{t}].}\leqno(21)
$$

\endproclaim

\noindent
{\bf Proof.} In the case of the brackets $\{\bar\xi_{\a\b},
\bar\xi_{\b\g}\}$,
$\a,\b,\g$ all different, it is a consequence of the formula
$$
{\th(tzw^{-1})\over{\th(t)\th(zw^{-1})}}
{\th(tt'w)\over{\th(tt')\th(w)}}
-{\th(t^{\prime -1}zw^{-1})\over{\th(t^{\prime -1})\th(zw^{-1})}}
{\th(tt'z)\over{\th(tt')\th(z)}} =
{\th(tz)\over{\th(t)\th(z)}}
{\th(t'w)\over{\th(t')\th(w)}};
$$
to show it, let $F(z,w,t,t')$ be the difference of both sides. We have
$F(qz,w,t,t')=t^{-1}F(z,w,t,t')$; moreover $F$ has no poles for $z\to 0$
or $z\to w$; since $t\notin q^{\ZZ}$, this shows $F=0$.

In the case of the brackets $\{\bar\xi_{\a\a},\bar\xi_{\a\b}\}$, $\a\ne
\b$,
it follows from
$$
\eqalign{
{1\over\th'(1)}t{d\over dt}[{\th(tw)\over{\th(t)\th(w)}}]
+{1\over\th'(1)}{\dot\th\over\th}(z){\th(tw)\over{\th(t)\th(w)}}
= &
-{\th(tzw^{-1})\over{\th(t)\th(zw^{-1})}}
{\th(tz)\over{\th(t)\th(z)}} \cr
&
+{1\over\th'(1)}{\dot\th\over\th}(zw^{-1}){\th(tw)\over{\th(t)\th(w)}};
}$$
this equality is proven as follows; let $F(z,w)$ be the difference of
the two sides, then $F(qz,w)=F(z,w)$ and $F(z,w)$ has  no poles for
$z\to w $ or $z\to 0$, which shows that $F(z,w)$ does not depend on $z$;
pose $F(z,w)=\varphi(w)$, then $\varphi(qw)=t\varphi(w)$ (this follows
from
${\dot\th\over\th}(qz)=-1+{\dot\th\over\th}(z)$, obtained by derivation
of the functional equation in $\th$); and $\varphi(z)$ has no poles
either, so $F(z,w)=0$.

In the case of the brackets $\{\bar\xi_{\a\b},\bar\xi_{\b\a}\}$, $\a\ne
\b$, it follows from the fact, that if we pose
$$
F(z,w)= {\th(t^{-1}z)\over{\th(t^{-1})\th(z)}}{\th(tw)\over{\th(t)\th(w)}}
+{\th(t^{-1}zw^{-1})\over{\th(t^{-1})\th(zw^{-1})}}{1\over\th'(1)}
[{\dot\th\over\th}(z)-{\dot\th\over\th}(w)],
$$
we have $F(z,w)=F(z\zeta,w\zeta)$ for any $\zeta\in\CC^{\times}$.
Indeed, $F(qz,w)=tF(z,w)-{1\over\th'(1)}t
{\th(t^{-1}zw^{-1})\over{\th(t^{-1})\th(zw^{-1})}}$,
so with $\varphi(z,w)=F(z,w)-F(z\z,w\z)$, we have
$\varphi(qz,w)=t\varphi(z,w)$; as $F(z,w)$ has no poles in $z$ (or in
$w$), $\varphi$ has no poles either, and so it vanishes. So $F$ is only
a function of $zw^{-1}$, that we can evaluate when $w\to 1$; this
evaluation gives the matrix elements of $\rho$.

The brackets $\{\bar\xi_{\a\a},\bar\xi_{\b\b}\}$ are all zero, and the
$[r,\bar\xi\otimes 1 + 1\otimes\bar\xi]_{\a\a}^{\b\b}$ also; finally, the
brackets $\{\bar\xi_{\a\b},\bar\xi_{\g\d}\}$ ($\a,\b,\g,\d$ all
different) are
all zero, as well as the matrix elements
$[r,\bar\xi\otimes 1 + 1\otimes\bar\xi]_{\a\g}^{\b\d}$.  \hfill $\bull$

Now, after the reduction by $T$, the $\tr\bar\xi(z)^{s}$ will be
in involution.

Let us give now the explicit form of the Hamiltonians generated by
$\tr\bar\xi(z)^{2}$. We have
$$\eqalign{
\tr\bar\xi(z)^{2} & =\sum_{\a=1}^{n}(p_{\a}+\sum_{i=1}^{N}\eta_{\a\a}^{(i)}
{\dot\th\over\th}(zz_{i}^{-1}))^{2}
+2\sum_{1\le \a<\b\le n}
\bigg(
\sum_{i=1}^{N}\eta_{\a\b}^{(i)}{\th(t_{\a}t_{\b}^{-1}zz_{i}^{-1})\over
\th(t_{\a}t_{\b}^{-1})\th(zz_{i}^{-1})}
\bigg)
\cdot \cr
 &\cdot
\bigg(
\sum_{i=1}^{N}\eta_{\b\a}^{(i)}{\th(t_{\b}t_{\a}^{-1}zz_{i}^{-1})\over
\th(t_{\b}t_{\a}^{-1})\th(zz_{i}^{-1})}
\bigg);\cr}
$$
since
$$
\bar\xi(qz)=\Ad(t_{\a})\bar\xi(z)-(\sum_{i=1}^{N}\eta^{(i)})_{t},
$$
we have
$$
(\tr\bar\xi^{2})(qz)=(\tr\bar\xi^{2})(z)+\tr(\sum_{i=1}^{N}
\eta^{(i)})_{t}^{2}-2\sum_{i=1}^{N}{{\dot\theta}\over{\theta}}(zz_{i}^{-1})
\{
\sum_{\a=1}^{n}\eta_{\a\a}(\sum_{i}\eta_{\a\a}^{(i)})\},
$$
so that
$$
\tr\bar\xi(z)^{2}=H_{0}+\sum_{i=1}^{N}H_{i}
{{\dot\theta}\over{\theta}}(zz_{i}^{-1})+\sum_{i=1}^{N}\bigg(
{{\dot\theta}\over{\theta}}(zz_{i}^{-1})\bigg)^{2}
\{\sum_{\a=1}^{n}\eta_{\a\a}(\sum_{i}\eta_{\a\a}^{(i)})\};
$$
using
$$
\eqalign{
({\dot\th\over\th}(zz_{i}^{-1})-{\dot\th\over\th}(zz_{j}^{-1}))^{2}
 & =\wp(\ln (zz_{i}^{-1}))
+\wp(\ln (zz_{j}^{-1})) \cr
&
-2{\dot\th\over\th}(z_{i}z_{j}^{-1})
[{\dot\th\over\th}(zz_{i}^{-1})-{\dot\th\over\th}(zz_{j}^{-1})]
+[{\dot\th\over\th}(z_{i}z_{j}^{-1})]^{2},
}
$$
$$\eqalign{
{\th(tzz_{i}^{-1})\over{\th(t)\th(zz_{i}^{-1})}}
{\th(t^{-1}zz_{j}^{-1})\over{\th(t^{-1})\th(zz_{j}^{-1})}}
 & =[{\dot\th\over\th}(zz_{i}^{-1})-{\dot\th\over\th}(zz_{j}^{-1})]
{\th(t^{-1}z_{i}z_{j}^{-1})\over{\th(t^{-1})\th(z_{i}z_{j}^{-1})}} \cr
 & -[{\dot\th\over\th}(t^{-1})-{\dot\th\over\th}(t^{-1}z_{i}z_{j}^{-1})]
{\th(t^{-1}z_{i}z_{j}^{-1})\over{\th(t^{-1})\th(z_{i}z_{j}^{-1})}}
{\rm\  if\ } j\ne i, \cr
& =\wp(\ln (zz_{i}^{-1}))-\wp(\ln (tz_{i})) {\rm\ \  else }
\cr}$$
[we set $z=e^{\tau}$, so near $\tau=0$,
${\dot\th\over\th}(z)\sim{1\over\tau}$,
$\wp(\tau)\sim{1\over{\tau^{2}}}$], we find
$$\eqalign{
H_{i} & =2\sum_{\a=1}^{n}p_{\a}\eta_{\a\a}^{(i)}+2\sum_{\a=1}^{n}
\sum_{j\ne i}\eta_{\a\a}^{(i)}
\eta_{\a\a}^{(j)}
[{\dot\th\over\th}(z_{i}z_{j}^{-1})-{\dot\th\over\th}(z_{j}z_{i}^{-1})] \cr
 & +2\sum_{\a\ne \b}\sum_{j\ne i}\eta_{\a\b}^{(i)}\eta_{\b\a}^{(j)}
{\th(t_{\b}t_{\a}^{-1}z_{i}z_{j}^{-1})\over{\th(t_{\b}t_{\a}^{-1})}\th(z_{i}
z_{j}^{-1})} \cr
}\leqno(22)
$$
(a less symmetric form could be obtained using
${\dot\th\over\th}(z)+{\dot\th\over\th}(z^{-1})=1$, and the irrelevance
of combinations of the
$\sum_{i=1}^{N}\eta_{\a\a}^{(i)}$), and
$$
\eqalign{
H_{0} & =\sum_{\a=1}^{n}p_{\a}^{2}+\sum_{j<i} \eta_{\a\a}^{(i)}
\eta_{\a\a}^{(j)}[{\dot\th\over\th}(z_{i}z_{j}^{-1})]^{2}
-2\sum_{\a<\b}\sum_{i=1}^{N} \eta_{\a\b}^{(i)}
\eta_{\b\a}^{(i)}\wp(\ln(t_{\a}t_{\b}^{-1})) \cr
&
-2\sum_{\a<\b}\sum_{i\ne j}\eta_{\a\b}^{(i)}\eta_{\b\a}^{(j)}
[{\dot\th\over\th}(t_{\b}t_{\a}^{-1})-{\dot\th\over\th}
(t_{\b}t_{\a}^{-1}z_{i}z_{j}^{-1})]
{\th(t_{\b}t_{\a}^{-1}z_{i}z_{j}^{-1})\over{\th(t_{\b}t_{\a}^{-1})
\th(z_{i}z_{j}^{-1})}}.
\cr}\leqno(23)
$$
\section{Remark.}{} It is interesting to compare these results with
those of [15], [5], [19], [11]. The system considered in these papers is
connected with the case $N=1$. Also there should be some connection
between the $r$-matrix (19) and the ones from [5] and [19].

\section{5.}{Gaudin-Calogero system in the $sl_{2}$ case.}
Let us see now how to construct a quantization of the system of the last
section when $G=GL_{2}$. We will construct differential operators on
$$
\cM^{(0)}_{G}(X)=
T\rtimes S_{n}\setminus (\CC^{\times})^{n}\times
(G/B)^{N}/[(t_{\a},g_{i}B)
\sim (q^{a_{\a}}t_{\a},\diag(z_{i}^{a_{\a}})g_{i}B)],
$$
whose symbols will be the Hitchin's Hamiltonians, $\tr\bar\xi(z)^{s}$,
for $n=2$.

For this, we consider an integer $k$ and a system of dominant weights
$(\la_{i})_{i=1,\cdots,N}$, and the algebra
$\cA=\Diff((\CC^{\times})^{n},\cL_{k}^{\boxtimes n})\otimes
\otimes_{i=1}^{N}\Diff(G/B,\cL_{\la_{i}})$ [here
$\cL_{k}=\pi^{*}\cO(k(1))$, $\pi$ the projection $\CC^{\times}\to X$,
$\Diff(X,\cL)=H^{0}(X,\cL\otimes\cD_{X}\otimes\cL^{-1})$, for $X$ an
analytic variety and $\cL$ a line bundle on $X$]. Let $(t_{\a})_{1\le
\a\le n}$ be the coordinates on $(\CC^{\times})^{n}$, and
$\hat{p}_{\a}=t_{\a}{\pr\over{\pr t_{\a}}}+k{\dot\th\over\th}(t_{\a})$;
let again, $e_{\a\b}^{(i)}$
denote the action of $e_{\a\b}\in gl_{n}(\CC)$ on the $i$-th factor of
the second part of $\cA$. Consider now the matrix $L(z)\in
gl_{n}(\CC)\otimes \cA$, defined by
$$
\eqalign{
L(z)_{\a\b}=\sum_{i=1}^{N}  e_{\a\b}^{(i)}
{{\th(t_{\a}^{-1}t_{\b}z z_{i}^{-1})}\over {\th(t_{\a}^{-1}t_{\b})\th(z
z_{i}^{-1})}}  & {\rm \ if\  } \a\ne \b, \cr &
L(z)_{\a\a}={1\over\th'(1)}\hat{p}_{\a}
+\sum_{i=1}^{N}{1\over\th'(1)}{\dot\th\over\th}(zz_{i}^{-1})e_{\a\a}^{i}.
}\leqno(24)
$$

Let us perform now the reduction of $\cA$ w.r.t. $T$. It can be done as
follows: let $\cA[0]$ be the weight zero subalgebra of $\cA$, $\cA[0]=\{x\in
\cA| [h_{\a\a},x]=0, 1\le \a \le n\}$ and $\cA^{red}=\cA[0]/(h_{\a\a})_{1\le \a
\le n}$ (where  $(h_{\a\a})_{1\le \a\le n}$ is the left, or right ideal
generated by the $h_{\a\a}$ in $\cA[0]$). Then $\cA^{red}$ is the algebra of
globally defined differential operators on $(\CC^{\times})^{n}\times
[T\setminus(G/B)^{N}]$, twisted by the quotient of $\cL_{k}^{\boxtimes
n}\boxtimes\boxtimes_{i=1}^{N}
\cL_{\la_{i}}$.

{}From $\tr L(qz)^{2}=\tr L(z)^{2}+\tr(L(z)h+hL(z))+\tr(h^{2})$, we see
that $\tr L(z)^{s}$, $s=1,2$ define elements of $[\cA/\sum_{\a=1}^{2}
\cA h_{\a\a}]\otimes H^{0}(X,\cO(s\sum_{i=1}^{N}(z_{i})))$, which also
belong to
$\cA^{red}\otimes H^{0}(X,\cO(s\sum_{i=1}^{N}(z_{i})))$. Then

\proclaim{Proposition 5.1} The expansions of $\tr L(z)^{s}$, $s=1,2$,
along bases of $ H^{0}(X,\cO(s\break\sum_{i=1}^{N}(z_{i})))$, form a
commutative family in $\cA^{red}$. These operators are $S_{2}$-invariant
and invariant under the action of $\ZZ^{2}$ defined by
$(a_{\a})\cdot
(t_{\a},g_{i}B)=(q^{a_{\a}}t_{\a},\diag(z_{i}^{a_{\a}})g_{i}B)$, and hence
define operators on
$\cM^{(0)}_{GL_{2}}(X)$, twisted by the line bundle associated
with $(k,\la_{i})$. Their symbols coincide with Hitchin's Hamiltonians.
\endproclaim

\noindent
{\bf Proof.} If $w\in S_{2}$, then $w^{*}L(z)=L(z)$; let
$\vare_{\a}$ be the $\a$-th basis vector of $\ZZ^{2}$, then
$\vare_{\a}^{*}L(z)=\Ad(\diag(1,\cdots,z,\cdots,1))L(z)$ ($z$ in $\a$-th
position). The last statement follows from the fact that the symbol of
$\hat{p}_{\a}$ is $p_{\a}$, and the symbol of $e_{\a\b}^{(i)}$ is
$\eta_{\a\b}^{(i)}$. The first statement follows from a direct
computation, using the explicit form of the Hamiltonians:
$$
\eqalign{
\hat{H}_{i}
 & = \hat{p}h^{(i)}
+2\sum_{j\ne i}h^{(i)}h^{(j)}
[{\dot\th\over\th}(z_{i}z_{j}^{-1})-{\dot\th\over\th}(z_{j}z_{i}^{-1})]
\cr
 & +2\sum_{j\ne i} \big(e^{(i)}f^{(j)}
{\th(t^{2}z_{i}z_{j}^{-1})\over{\th(t^{2})}\th(z_{i}
z_{j}^{-1})}
+e^{(j)}f^{(i)}
{\th(t^{-2}z_{i}z_{j}^{-1})\over{\th(t^{-2})}\th(z_{i}
z_{j}^{-1})}\big)
}\leqno(25)
$$
and
$$
\eqalign{
\hat{H}_{0} & =\hat{p}^{2}+\sum_{j<i} h^{(i)}
h^{(j)}[{\dot\th\over\th}(z_{i}z_{j}^{-1})]^{2}
-2\sum_{i=1}^{N}e^{(i)}f^{(i)}\wp(\ln(t^{2})) \cr
&
-2\sum_{i\ne j}e^{(i)}f^{(j)}
[{\dot\th\over\th}(t^{2})-{\dot\th\over\th}
(t^{2}z_{i}z_{j}^{-1})]
{\th(t^{2}z_{i}z_{j}^{-1})\over{\th(t^{2})
\th(z_{i}z_{j}^{-1})}}.
\cr}\leqno(26)
$$
\hfill $\bull$

Equations (25) and (26) define differential operators acting on
$\CC^{\times}\times
[T\setminus (\CC P^{1})^{N}]$; $(t,t_{i})$
being the product coordinates on
$\CC^{\times}\times (\CC P^{1})^{N}$, we have
$\hat{p}=2t{\pr\over{\pr t}}+2k{\dot\th\over\th}(t^{2})$,
$h^{(i)}=2(t_{i}{\pr\over{\pr t_{i}}}+\la_{i})$,
$e^{(i)}=t_{i}^{2}{\pr\over{\pr t_{i}}}+2\la_{i}t_{i}$,
$f^{(i)}=-{\pr\over{\pr t_{i}}}$.

For $N=1$, this system is reduced to
$$
\hat{H}_{0}=\hat{p}^{2}-2e^{(1)}f^{(1)}\wp(\ln t^{2}) ,
\quad\hat{H}_{1}=e^{(1)}f^{(1)}.
$$

\section{Remarks.}{}
1. A natural module for $\cA^{red}$ is ${\rm Fun}(\CC^{\times})\otimes
(V_{\la_{1}}\otimes\cdots\otimes V_{\la_{N}})[0]$.
More precisely, we can pose the eigenvalue problem
$\hat{H}_{i}\psi=\mu_{i}\psi$, $\hat{H}_{0}\psi=\mu_{0}\psi$, $\psi$ a
function of $\CC^{\times}$, with values in
$\otimes_{i=1}^{N}V_{\la_{i}}$, whose component in
$(\otimes_{i=1}^{N}V_{\la_{i}})[\bar\la_{i}]$, $\psi_{\bar\la_{i}}(t)$,
satisfies $\psi_{\bar\la_{i}}(qt)=z_{1}^{\bar\la_{1}}\cdots
z_{N}^{\bar\la_{N}}z^{\ell}
\psi_{\bar\la_{i}}(t)$, for each system of weights $(\bar\la_{i})$,
$\ell$ being a fixed integer.
The space of such functions, with only poles at
$q^{\ZZ}$, is stable under the actions of $\hat{H}_{0}$ and the
$\hat{H}_{i}$.

2. Prop. 5.1 suggests that the operators constructed here coincide with
the result of the action of the center of the enveloping algebra at the
critical level, when $k=2$. Indeed in this case, after [22], the
quotient of $\cL_{k}^{\boxtimes 2}$ by $S_{2}$ coincides with
$(\det_{|\cM^{(0)}_{GL_{2}}(X)})^{-2}$, on which this center should act.

After obtaining the main results of this paper, we learnt about the
paper of N. Nekrasov [16], where Hitchin systems for degenerate curves
are described as many-body problems.

\vskip 1truecm
\noindent
{\bf References}
\bigskip
\item{[1]} M.R. Adams, J. Harnad, E. Previato, {\sl Isospectral
Hamiltonian flows in finite and infinite dimensions II. Integration of
flows,} Commun. Math. Phys. 134 (1990), 555-85.
\medskip
\item{[2]} M. Adler, {\sl On a trace functional for formal
pseudo-differential operators and the symplectic structure for the KdV
type equations,} Invent. Math. 50 (1979), 219-48.
\medskip
\item{[3]} A. Beauville, {\sl Jacobiennes des courbes spectrales et
syst\`emes compl\`etement int\'e-grables,} Acta Math., 169 (1990),
211-35.
\medskip
\item{[4]} A.A. Beilinson, V.G. Drinfeld, {\sl Quantization of
Hitchin's fibration and Langlands program,} preprint.
\medskip
\item{[5]} H.W. Braden, T. Suzuki, {\sl $R$-matrices for Elliptic
Calogero-Moser Models}, Lett. Math. Phys. 30, 147-59 (1994).
\medskip
\item{[6]} B.L. Feigin, E.V. Frenkel, N. Reshetikhin, {\sl Gaudin
model, Bethe ansatz and correlation functions at the critical level,}
Commun. Math. Phys. 166 (1), 27-62 (1995).
\medskip
\item{[7]} G. Felder, C. Wieczerkowski, {\sl Conformal field theory on
elliptic curves and Knizhnik-Zamolodchikov-Bernard equations,}
hep-th/9411004.
\medskip
\item{[8]} R. Garnier, Rend. Circ. Mat. Palermo 43, 155-91 (1919).
\medskip
\item{[9]} M. Gaudin, Jour. Physique 37 (1976), 1087-1098.
\medskip
\item{[10]} R. Goodman, N.R. Wallach, {\sl Higher-order Sugawara operators
for affine Lie algebras,} Trans. AMS, 315:1 (1989), 1-55.
\medskip
\item{[11]} A.S. Gorsky, N.A. Nekrasov,{\sl Elliptic Calogero-Moser
system from two-dimensio-nal current algebra,} hep-th/9401021.
\medskip
\item{[12]} N. Hayashi, {\sl Sugawara operators and Kac-Kazhdan
conjecture,} Invent. Math. 54 (1988), 13-52.
\medskip
\item{[13]} N. Hitchin, {\sl Stable bundles and integrable systems,} Duke
Math. Jour., 54 (1), 91-114 (1987).
\medskip
\item{[14]} B. Kostant, {\sl The solution to a generalized Toda lattice
and representation theory,} Adv. Math. 34 (1980), 13-53.
\medskip
\item{[15]} I.M. Krichever, {\sl Elliptic solutions of the
Kadomtsev-Petviashvili equation and integrable systems of particles,}
Funct. An. Appl., 14 (1), 282-90 (1990).
\medskip
\item{[16]} N. Nekrasov, {\sl Holomorphic bundles and many-body systems,}
PUPT-1534, ITEP-N95/1, hep-th/9503157.
\medskip
\item{[17]} A.G. Reyman, {\sl Quantum tops,} Int. J. Mod. Phys. B,
7:20-21 (1993), 3707-13.
\medskip
\item{[18]} M.A. Semenov-Tian-Shansky, {\sl D. Sci. thesis,} LOMI,
Leningrad (1985).
\medskip
\item{[19]} E.K. Sklyanin, {\sl Dynamical $r$-matrices for the elliptic
Calogero-Moser system,} LPTHE 93-42, hep-th/9308060.
\medskip
\item{[20]} T.A. Springer, {\sl Trigonometric sums, Green functions of
finite groups and representations of Weyl groups,}
Inv. Math. 36, 173-207 (1976).
\medskip
\item{[21]} W. Symes, {\sl Systems of Toda type, inverse spectral
problems and representation theory,} Invent. Math. 59 (1990),
195-338.
\medskip
\item{[22]} L.W. Tu, {\sl Semistable bundles over an Elliptic Curve,}
Adv. Math. 98, 1-26 (1993).
\medskip
\medskip\medskip
\section{}{}

B.E., V.R.: Centre de Math\'{e}matiques, URA 169
du CNRS, Ecole Polytechnique, 91128 Palaiseau, France

V.R.: ITEP, Bol. Cheremushkinskaya, 25, 117259,
Moscow, Russia.
\bye